\documentclass[11pt,reqno]{amsart}
\usepackage[T1]{fontenc}
\usepackage{lmodern,microtype}
\usepackage[a4paper,margin=29mm]{geometry}
\usepackage{amssymb,mathtools}
\usepackage{xcolor}
\usepackage[colorlinks=true,linkcolor=blue!45!black,citecolor=blue!45!black,
  urlcolor=blue!45!black]{hyperref}
\hypersetup{pdftitle={On the sharpness of a Liouville exponent in a competitive elliptic
Schr\"odinger system},pdfauthor={}}
\newtheorem{theorem}{Theorem}
\newtheorem{lemma}{Lemma}[section]
\theoremstyle{remark}
\newtheorem*{remark}{Remark}
\newcommand{\R}{\mathbb R}
\newcommand{\T}{\mathbb T}
\newcommand{\cH}{\mathcal H}
\newcommand{\avg}[1]{\overline{#1}}
\newcommand{\dt}{\,dt}

\newcommand{\ess}{\mathrm{ess}}
\numberwithin{equation}{section}
\allowdisplaybreaks[1]
\title[On the sharpness of a Liouville exponent]
{On the sharpness of a Liouville exponent in a competitive elliptic
Schr\"odinger system}

  
\date{\today}

\author{Fengmiao Bian}
\address[Fengmiao Bian]{%
  Department of Mathematics,
  The Hong Kong University of Science and Technology,
  Clear Water Bay, Hong Kong, China}
\email{mafmbian@ust.hk}

\author{Hui Yang}
\address[Hui Yang]{%
  School of Mathematical Sciences, 
  Shanghai Jiao Tong University, Shanghai 200240, China}
\email{hui-yang@sjtu.edu.cn}

\begin{document}
\begin{abstract}
We identify \(p_*(n)=(n-1)/(n-3)\), \(n\ge5\), as the sharp exponent
for Liouville theorems concerning a two-component elliptic Schr\"odinger system with
self-focusing interactions and repulsive intercomponent coupling. We consider the strictly copositive regime, in which the repulsion
strength is smaller than the geometric mean of the two self-interaction strengths. Quittner and Souplet (Comm. Math. Phys. 2012: 1--19) established a Liouville-type nonexistence result in the range \(1<p<p_*(n)\). For every \(p_*(n)<p<(n+2)/(n-2)\), we construct bounded, nontrivial, nonnegative entire solutions whenever the repulsion strength approaches the copositivity threshold from below. 
%These solutions are periodic in one coordinate, decay in the remaining coordinates, and have finite Dirichlet energy on each period cell but infinite energy on the whole space. 
The construction combines a compact variational functional, a periodic spectral gap, and uniform exterior estimates that remove a truncation. Our result gives a negative answer to the open question posed by Quittner and Souplet (Comm. Math. Phys. 2012: 1--19) of whether strict copositivity implies nonexistence of nontrivial nonnegative entire solutions throughout the Sobolev subcritical range when \(n\ge5\).
\end{abstract}
\maketitle

\section{Introduction and main results} 
In this paper, we study nonnegative entire solutions of the coupled elliptic Schr\"odinger system
\begin{equation}\label{eq:system}
\begin{cases}
-\Delta u=\mu_1u^{2q+1}+\beta u^qv^{q+1},\\
-\Delta v=\mu_2v^{2q+1}+\beta v^qu^{q+1} 
\end{cases}
\quad\text{in }\R^n, 
\end{equation}
where \(q>0\), \(\mu_1,\mu_2>0\), and \( -\sqrt{\mu_1\mu_2} < \beta <0\).
Coupled nonlinear Schr\"odinger equations describe interacting optical
waves in nonlinear media and multicomponent Bose--Einstein condensates \cite{Frantzeskakis,Kivshar}; their standing-wave profiles satisfy elliptic systems, usually with additional linear frequency and trapping-potential terms \cite{DWW,LW,Sirakov}. The homogeneous system above also arises as a limiting equation under rescaling near concentration points, where these lower-order
terms disappear \cite{DWW,PQS}.
Here \(u\) and \(v\) represent the amplitudes of the two components,
while \(\mu_1\) and \(\mu_2\) measure their self-interaction strengths
and \(\beta\) measures the intercomponent coupling. With the sign convention used here, positive self-interaction coefficients correspond to self-focusing in optics and attractive intraspecies interactions in condensates, whereas negative coupling describes repulsion between the components. The exponent \(p=2q+1\) specifies the degree of the nonlinear response
in the wave amplitudes. The case \(q=1\), or \(p=3\), recovers the cubic nonlinearities
associated with Kerr optics and the Gross--Pitaevskii description
of dilute condensates. General powers provide a phenomenological extension allowing
different dependences of the effective response on intensity or density \cite{Frantzeskakis,Kivshar}. 
The existence and multiplicity of solutions, as well as uniform
regularity estimates, have been extensively studied for the coupled elliptic Schr\"odinger system and related models; see, for example,
\cite{AC07,BDW10,CZ12,CZ15,DWW,GLW,LW,NTTV10,PW13,Sirakov,SZ15,TTVW,TV09,WW08}.
A priori estimates and Liouville-type theorems have also been widely investigated; see \cite{BDW10,DWW,QS12,QS25,Souplet,TTVW} and the references therein.

Throughout the paper, nontrivial means \((u,v)\not\equiv(0,0)\). We write
\[
p=2q+1,\qquad p_*(n)=\frac{n-1}{n-3}, \qquad p_S(n)=\frac{n+2}{n-2}. 
\]
The coupling matrix \(B=\left(\begin{smallmatrix}\mu_1&\beta\\
\beta&\mu_2\end{smallmatrix}\right)\) is strictly copositive precisely when \(\beta>-\sqrt{\mu_1\mu_2}\). Strict copositivity plays a significant role in quadratic programming \cite{Jacobson}. In the competitive range considered here, the matrix is then even
positive definite. Nevertheless, a component of the right-hand side of \eqref{eq:system} may change sign. For the scalar equation \(-\Delta w=w^p\) on \(\R^n\), \(n\ge3\),
Gidas and Spruck \cite[Theorem~1.1]{GS} proved that \(1<p<p_S(n)\)
excludes every nontrivial nonnegative classical solution. This gives the corresponding conclusion for
\eqref{eq:system} at \(\beta=0\) and for semitrivial solutions at any
\(\beta\). For cooperative coupling \(\beta>0\), the Liouville theorem
holds in the same full subcritical range by the method of moving spheres or moving planes for cooperative systems \cite{ChenLi,GuoLiu,RZ,Zou}.

For competing components, several distinctions in the earlier results
are important here.
In the cubic case \(q=1\), Dancer, Wei and Weth
\cite{DWW} proved triviality of nonnegative entire
solutions of \eqref{eq:system} in dimensions \(n\le3\) when
\(\beta>-\sqrt{\mu_1\mu_2}\).
Tavares, Terracini, Verzini and Weth \cite{TTVW} treated general
homogeneous powers as well. For two components and \(n\ge3\),
their result gives triviality under strict copositivity for
\(1<p\le n/(n-2)\). For \(\beta>-\sqrt{\mu_1\mu_2}\), Quittner and Souplet
\cite{QS12} improved the results in \cite{TTVW} and established nonexistence for nonnegative classical entire
solutions of \eqref{eq:system} when 
\[
\begin{cases}
1<p<p_S(n), & n=3,4\\
1<p<p_*(n), & n\ge5. 
\end{cases} 
\] 
Their method uses a Rellich--Pohozaev identity and interpolation on
spheres. For a subsequent general gradient-system formulation, see
\cite{QS25}. The multicomponent result of Souplet \cite{Souplet} concerns \(n\ge3\), at least
three components, a strictly copositive matrix and \(1<p<n/(n-2)\).

Our first result establishes the sharpness of the exponent
\(p_*(n)=(n-1)/(n-3)\), \(n\ge5\), in the coefficient-uniform
Liouville theorem of Quittner and Souplet \cite{QS12}.
Their theorem excludes nontrivial nonnegative entire solutions
for \(1<p<p_*(n)\) under strict copositivity.
For every \(p_*(n)<p<(n+2)/(n-2)\), we construct bounded nonnegative
entire solutions with both components nontrivial, for equal
positive self-interaction coefficients and competitive coupling
sufficiently close to the copositivity boundary. These solutions provide the complementary existence result and
show that the nonexistence interval in \cite{QS12} cannot be
enlarged uniformly over the admissible coefficients, even when
attention is restricted to bounded solutions.
Thus \(p_*(n)\) is the optimal upper threshold for this uniform
Liouville property.

\begin{theorem}[Existence near the copositivity boundary]\label{thm:exist}
Let \(n\ge5\) and \(\mu_1=\mu_2=\mu>0\). Suppose that 
\begin{equation}\label{eq:gap}
 p_*(n)<p<p_S(n), \qquad \text{equivalently}\qquad \frac1{n-3}<q<\frac2{n-2}. 
\end{equation}
Then there exists \(\eta=\eta(n,p)\in(0,1)\) such that for every
\(\beta\in(-\mu,-\mu+\mu\eta)\), \eqref{eq:system} admits a bounded nonnegative
\(C^2\) solution with both components nonzero and \(u\not\equiv v\).
Writing \(x=(y,t)\in\R^{n-1}\times\R\), the components can be chosen
radial in \(y\), \(2\pi\)-periodic in \(t\), and related by
\(v(y,t)=u(y,t+\pi)\). For some \(\gamma>1/q\),
\[
 u(y,t)+v(y,t)\le C(1+|y|)^{-\gamma}.
\]
The Dirichlet energy is finite on one period cell and infinite on
\(\R^n\).
\end{theorem}

Theorem~\ref{thm:exist} gives a
negative answer to the open question raised by Quittner and
Souplet \cite{QS12} of whether strict copositivity
ensures the Liouville property for nonnegative entire solutions
throughout the Sobolev subcritical range in dimensions \(n\ge5\). 
We now introduce the idea behind the construction in Theorem~\ref{thm:exist}. We work on \(\R^{n-1}\times\T\) and impose a half-period translation
that exchanges the two components. Their period averages then agree,
so the functional \(\int(|u|^{q+1}-|v|^{q+1})^2\) only measures a
difference between the components. Radial symmetry in the unbounded
variables and cancellation of the common mean make this functional
compact on the Dirichlet unit ball. We add a truncated positive
pure-power term to move the coupling into the strictly copositive
range. The estimates needed to remove the truncation are uniform in
its parameter. The sum of the components is positive and satisfies a
uniform Harnack inequality. Their difference is antiperiodic and
therefore has a spectral gap in the periodic variable, which gives
exponential decay of its squared period average. The common mean
then satisfies a scalar differential inequality with an exponentially
decaying source. An exterior comparison, justified by Hardy's
inequality, gives decay faster than \(|y|^{-1/q}\). This decay removes
the truncation and yields a solution of the original system after an
amplitude rescaling.

The preceding construction produces nontrivial entire solutions
when the competitive coupling approaches the copositivity boundary.
We next examine the regime of weak competition, near the decoupled
system \(\beta=0\), where the scalar Liouville theorem excludes
nontrivial nonnegative entire solutions throughout the Sobolev
subcritical range.
Our second result shows that this nonexistence property persists
under sufficiently small competitive perturbations, for arbitrary
positive self-interaction coefficients.
Moreover, the coupling neighborhood can be chosen uniformly for
\(1<p \le\bar p< p_S(n)\), without any boundedness or decay assumptions
on the solutions.
Thus, above \(p_*(n)\), the existence question depends essentially on
the strength of the competition, even though a Liouville theorem
uniform over all admissible coefficients is no longer available.

\begin{theorem}[Liouville theorem for weak competition]\label{thm:weak}
Let \(n\ge5\), \(\mu_1,\mu_2>0\), and fix $1<\bar p<p_S(n)$ (equivalently, \(0<\bar q<2/(n-2)\)).
Then there exists \(\varepsilon_*=\varepsilon_*(n,\bar p,\mu_1,\mu_2) \in (0, \sqrt{\mu_1\mu_2}) \) such that for
every \(1<p\le\bar p\) and every \(-\varepsilon_*<\beta< 0\),
the only nonnegative \(C^2(\R^n)\) solution of \eqref{eq:system} is \(u=v=0\).
\end{theorem}

Theorem~\ref{thm:weak} establishes the stability of the scalar
Liouville property under weak competitive coupling, for arbitrary
positive self-interaction coefficients.
The smallness threshold is uniform for
\(1<p\le\bar p<p_S(n)\), including the limit \(p \downarrow1\),
and no a priori boundedness or decay assumptions are imposed.
The proof combines homogeneity with the doubling--rescaling
method of Pol\'a\v{c}ik, Quittner and Souplet \cite{PQS}:
a hypothetical sequence of nontrivial solutions with coupling
tending to zero yields a bounded, normalized solution of a
decoupled limiting system.
The contradiction follows from the scalar subcritical Liouville
theorem.

The remainder of the paper is organized as follows.
Section~2 introduces the symmetric variational framework on the
cylinder and establishes the compactness of the nonlinear functional.
In Section~3, we construct solutions of a truncated problem and
derive uniform decay estimates that allow us to remove the
truncation, completing the proof of Theorem~\ref{thm:exist}.
Section~4 proves Theorem~\ref{thm:weak} by a doubling--rescaling
and compactness argument.

\section{A compact functional on the cylinder}
We first prove Theorem~\ref{thm:exist} for \(\mu=1\); a common
amplitude rescaling at the end of the construction gives every
\(\mu>0\). Let
\[ a=q+1 \quad \text{and} \quad m=2a=p+1.\]
Set \(k=n-1\ge4\) and
\(\Omega=\R^k\times\T\), where \(\T=\R/(2\pi\mathbb Z)\).
All integrals without an indicated domain in this and the next
section are over \(\Omega\), with measure \(dy\,dt\).
Let \(\dot H^1(\Omega)\) be the completion of smooth periodic
functions of compact support in \(y\) under the Dirichlet norm.
Let \(\cH\subset\dot H^1(\Omega)^2\) be the closed subspace of
radial-in-\(y\) pairs \(U=(u,v)\) with \(v(y,t)=u(y,t+\pi)\).
Define
\begin{equation}\label{eq:EG}
 E(U)=\int (|\nabla u|^2+|\nabla v|^2),\qquad
 G(U)=\int (|u|^a-|v|^a)^2 .
\end{equation}
The absolute values allow unrestricted signed variations.

For \(z\in\dot H^1(\Omega)\) put
\(\avg z(y)=(2\pi)^{-1}\int_{\T}z(y,t)\dt\).
By Jensen's inequality and periodic Poincar\'e,
\(\avg z\in\dot H^1(\R^k)\) and \(z-\avg z\in H^1(\Omega)\),
with both norms bounded by a constant times \(\|\nabla z\|_2\).
In particular \(\avg z\in L^{2k/(k-2)}(\R^k)\).
These statements follow first for test functions and then by
completion; they also identify each element with a local \(H^1\)
function. The \(k\)-dimensional Hardy inequality, integrated in \(t\),
gives
\begin{equation}\label{eq:hardy}
 H_k\int |y|^{-2}|z|^2\le\int|\nabla_y z|^2,
 \qquad H_k=(k-2)^2/4.
\end{equation}

\begin{lemma}\label{lem:compact}
Under \eqref{eq:gap}, \(G\) is finite, bounded on the unit ball of
\(\cH\), and weakly sequentially continuous on bounded sets of \(\cH\).
\end{lemma}
\begin{proof}
It suffices to prove local compactness and uniform smallness of the
tails of the integrand. Let \(E(U)\le1\). The two means agree by the
half-period shift. Write \(u=w+u'\), \(v=w+v'\), and \(D=u-v=u'-v'\).
For \(r=|y|>0\), the radial estimate
\begin{equation}\label{eq:radial}
 |w(r)|\le C r^{-(k-2)/2}\|\nabla_y w\|_{L^2(\R^k)}
\end{equation}
follows by integrating the radial derivative from \(r\) to infinity
and applying Cauchy--Schwarz with weight \(\rho^{k-1}\). The limit at
infinity is zero since \(w\in L^{2k/(k-2)}\).

For \(z=u'\) or \(v'\), set \(S_j=\{j<|y|<j+1\}\times\T\) and
\(e_j=\int_{S_j}(|z|^2+|\nabla z|^2)\). Then \(\sum_{j\ge1}e_j\le C\).
To make the weights explicit, write \(z(y,t)=Z(r,t)\) and
\(c_k=|\mathbb S^{k-1}|\). On \(Q_j=(j,j+1)\times\T\),
\[
 e_j=c_k\int_{Q_j}
 (|Z|^2+|\partial_rZ|^2+|\partial_tZ|^2)r^{k-1}\,dr\,dt,
 \qquad
 \|Z\|_{H^1(Q_j)}^2\le c_k^{-1}j^{-(k-1)}e_j.
\]
Translation identifies every \(Q_j\) with \((0,1)\times\T\).
The two-dimensional Sobolev embedding into \(L^m\), for this fixed
finite \(m\), thus has a constant independent of \(j\). Since
\(j\le r\le j+1\le2j\),
\begin{equation}\label{eq:strip}
 \begin{split}
 \int_{S_j}|z|^m
 &\le c_k(2j)^{k-1}\|Z\|_{L^m(Q_j)}^m\\
 &\le Cj^{k-1}\bigl(j^{-(k-1)}e_j\bigr)^{m/2}
 =Cj^{-\alpha}e_j^{m/2},
 \qquad \alpha=(k-1)(m/2-1)>0.
 \end{split}
\end{equation}
For an integer \(N\ge1\), the inequality
\(\sum e_j^{m/2}\le(\sum e_j)^{m/2}\) gives
\begin{equation}\label{eq:strip-tail}
 \int_{|y|>N}|z|^m
 \le CN^{-\alpha}\left(\sum_{j\ge N}e_j\right)^{m/2}
 \le CN^{-\alpha}.
\end{equation}
Both \(u'\) and \(v'\) have these bounds. Since \(D=u'-v'\),
the same bound holds for \(D\).
The mean value inequality for \(s\mapsto |s|^a\) yields
\[
 (|u|^a-|v|^a)^2
 \le C\bigl(|w|^{2q}+|u'|^{2q}+|v'|^{2q}\bigr)|D|^2.
\]
The first term has tail at most \(CR^{-(k-2)q}\|D\|_2^2\).
For either fluctuation \(z\), H\"older's inequality gives
\[
 \int_{|y|>N}|z|^{2q}|D|^2
 \le
 \left(\int_{|y|>N}|z|^m\right)^{2q/m}
 \left(\int_{|y|>N}|D|^m\right)^{2/m}
 \le CN^{-\alpha},
\]
where \(2q+2=m\). Thus the tail of the integrand of \(G\) is
uniformly bounded by \(C(N^{-(k-2)q}+N^{-\alpha})\).
On a bounded cylinder the local \(H^1\) bound gives a uniform
\(L^m\) bound, since \(m<2n/(n-2)\). This also proves boundedness
of \(G\) on the unit ball.

If \(U_j\rightharpoonup U\) in \(\cH\), the continuous local
\(H^1\) embedding and Rellich compactness give strong local
\(L^m\) convergence. The limit is unique, so this conclusion holds
for the full sequence. In particular
\[
 \big\||u_j|^a-|u|^a\big\|_{L^2(Q)}
 \le C\bigl(\|u_j\|_{L^{2a}(Q)}^{a-1}
             +\|u\|_{L^{2a}(Q)}^{a-1}\bigr)
             \|u_j-u\|_{L^{2a}(Q)}\longrightarrow0
\]
on each bounded cylinder \(Q\), and the same holds for \(v_j\).
It follows that the integral defining \(G\) converges on \(Q\).
Combining this with the uniform tails proves \(G(U_j)\to G(U)\).
Homogeneity gives the assertions on every bounded set.
\end{proof}

\section{Construction and removal of a truncation}
\subsection*{The maximization problem}
Choose \(\gamma,\sigma\) so that
\begin{equation}\label{eq:choice}
 1/q<\gamma<k-2,\qquad 0<\sigma<\gamma(p-1)-2.
\end{equation}
The lower bound on \(q\) in \eqref{eq:gap} is exactly what makes
this possible. The upper bound gives \(\gamma>(k-2)/2\). Set
\[
 \omega(y)=(1+|y|^2)^{-(2+\sigma)/2}.
\]
Fix a smooth nonnegative pair \(U_0\in\cH\) with \(E(U_0)=1\) and
\(g_0=G(U_0)>0\). For instance take a common compactly supported
radial bump times \(1+\tfrac12\cos t\) and \(1-\tfrac12\cos t\), and
normalize. Fix \(e>0\) so small that \(e/(2H_k)\le g_0/(2m)\).
For \(0<\delta\le1\), define
\[
 h_\delta(y,s)=s\min\{\delta|s|^{p-1},e\omega(y)\},
 \qquad H_\delta(y,s)=\int_0^s h_\delta(y,\tau)\,d\tau.
\]
In particular \(0\le H_\delta(y,s)\le e\omega(y)s^2/2\).
Maximize
\begin{equation}\label{eq:functional}
 J_\delta(U)=\frac{1}{m}G(U)+\int [H_\delta(y,u)+H_\delta(y,v)]
 \qquad\text{over } U\in\cH,\quad E(U)\le1.
\end{equation}
Both terms are weakly sequentially continuous on this ball. For the
second term, strong \(L^2\) convergence on compact sets suffices,
since \(|h_\delta(y,s)|\le e|s|\); its tail is uniformly \(O(R^{-\sigma})\)
by \eqref{eq:hardy} and \(\sup_{|y|>R}|y|^2\omega(y)\le R^{-\sigma}\).
The first term is covered by Lemma~\ref{lem:compact}. A maximizer
therefore exists. Its value is at least \(g_0/m\), while the integral
of the two \(H_\delta\)'s is at most \(e/(2H_k)\). Hence
\begin{equation}\label{eq:mass}
 G(U)\ge g_0/2.
\end{equation}
It follows that \(E(U)=1\), since otherwise radial dilation in
amplitude up to the unit sphere would strictly increase \(G\) and
would not decrease the other term. Replacing both components by
their absolute values preserves \(E,G,J_\delta\) and the
symmetries. We take a nonnegative maximizer.

\subsection*{The multiplier and the sum equation}
Write \(U=(u,v)\) for this maximizer. For a smooth compactly
supported symmetric test pair \(\Phi=(\phi,\psi)\), put
\[
 I_\Phi=\int(\nabla u\cdot\nabla\phi+\nabla v\cdot\nabla\psi),
 \qquad c(s)=E(U+s\Phi)^{-1/2}.
\]
Then \(c(0)=1\) and \(c'(0)=-I_\Phi\). The curve
\(c(s)(U+s\Phi)\) lies on the constraint sphere for small \(s\).
Compactly supported derivatives of \(G\) are justified by local
\(L^m\) integrability. Using \(G(cV)=c^mG(V)\) handles the
amplitude derivative without assuming global differentiability of
\(G\) on \(\cH\). For the \(H_\delta\) terms, the derivative along
this curve is dominated, for small \(s\), by a constant times
\(e\omega(|U|^2+|\Phi|^2)\), which is integrable by
\eqref{eq:hardy}. Write
\[
 F_1=u^q(u^a-v^a)+h_\delta(y,u),\qquad
 F_2=v^q(v^a-u^a)+h_\delta(y,v).
\]
Differentiating at the maximum gives the exact identity
\[
 0=\int(F_1\phi+F_2\psi)
   -I_\Phi\left\{G(U)+\int[u h_\delta(y,u)+v h_\delta(y,v)]\right\}.
\]
The expression in braces is positive by \eqref{eq:mass}.
Its reciprocal, denoted by \(\lambda_\delta\), gives the equations
\begin{equation}\label{eq:pen}
\begin{aligned}
 -\Delta u&=\lambda_\delta\big[u^q(u^a-v^a)+h_\delta(y,u)\big],\\
 -\Delta v&=\lambda_\delta\big[v^q(v^a-u^a)+h_\delta(y,v)\big],\\
 \lambda_\delta^{-1}&=G(U)+\int [u h_\delta(y,u)+v h_\delta(y,v)] .
\end{aligned}
\end{equation}
Indeed the derivative of \(G/m\) in the first component is
\(u^q(u^a-v^a)\), since \(m=2a\).
Here is the extension from symmetric tests to arbitrary tests.
For \(Q\in O(k)\) define
\[
 (R_Q\Phi)(y,t)=\Phi(Qy,t),\qquad
 S(\phi,\psi)(y,t)=(\psi(y,t+\pi),\phi(y,t+\pi)).
\]
With normalized Haar measure \(dQ\) on \(O(k)\), set
\begin{equation}\label{eq:group-average}
 \mathcal A\Phi=\frac12\int_{O(k)}(R_Q\Phi+S R_Q\Phi)\,dQ.
\end{equation}
The operators \(R_Q\) and \(S\) commute, and \(S^2\) is the identity.
Consequently \(\mathcal A\Phi\) is a smooth compactly supported
element of \(\cH\). Define the linear residual
\[
 \mathcal L(\Phi)=I_\Phi-\lambda_\delta\int(F_1\phi+F_2\psi).
\]
Changes of variables show that
\(\mathcal L(R_Q\Phi)=\mathcal L(S\Phi)=\mathcal L(\Phi)\).
For this assertion, use the radiality of \(u,v,\omega\) and the
identities \(u(y,t+\pi)=v(y,t)\) and
\(F_1(y,t+\pi)=F_2(y,t)\). Thus
\(\mathcal L(\Phi)=\mathcal L(\mathcal A\Phi)=0\), proving
\eqref{eq:pen} against every smooth compactly supported periodic
test pair. By \eqref{eq:mass},
\begin{equation}\label{eq:lambda}
 0<\lambda_\delta\le 2/g_0 .
\end{equation}

In the remainder of the construction all upper bounds are uniform
for \(0<\delta\le1\). Put \(W=u+v\) and \(D=u-v\). Addition in
\eqref{eq:pen} gives
\begin{equation}\label{eq:sum}
 -\Delta W=\lambda_\delta\big[(u^q-v^q)(u^a-v^a)
                  +h_\delta(y,u)+h_\delta(y,v)\big].
\end{equation}
In particular \(0\le-\Delta W\le C W^p\). We give the uniform local
estimates in detail. Fix \(0<\rho<\pi/8\). Balls of radius \(4\rho\)
on \(\Omega\) lift isometrically to Euclidean \(n\)-balls.
Set \(s_k=2k/(k-2)\) and \(2_n^*=2n/(n-2)\).
The mean \(w\) and the fluctuations in Lemma~\ref{lem:compact} satisfy
\[
 \|w\|_{L^{s_k}(\R^k)}
 +\|u'\|_{L^2(\Omega)}+\|v'\|_{L^2(\Omega)}\le C.
\]
For every \(x_0=(y_0,t_0)\), H\"older's inequality gives
\[
 \int_{B_{4\rho}(x_0)}|w(y)|^2\,dy\,dt
 \le 8\rho\,|B^k_{4\rho}|^{1-2/s_k}
             \|w\|_{L^{s_k}(\R^k)}^2\le C.
\]
Together with \(E(U)=1\), this proves
\(\|W\|_{H^1(B_{4\rho}(x_0))}\le C\), with \(C\) independent of
\(x_0\) and \(\delta\). The Euclidean Sobolev inequality on these
fixed-radius balls therefore gives
\(\|W\|_{L^{2_n^*}(B_{4\rho}(x_0))}\le C\).
Write \(-\Delta W=BW\), defining \(B=0\) on \(\{W=0\}\).
Then \(0\le B\le CW^{p-1}\), and
\[
 \ell=\frac{2_n^*}{p-1}>\frac n2,\qquad
 \|B\|_{L^\ell(B_{4\rho}(x_0))}
 \le C\|W\|_{L^{2_n^*}(B_{4\rho}(x_0))}^{p-1}\le C.
\]
The interior local boundedness estimate for this linear equation
gives
\[
 \ess\sup_{B_{2\rho}(x_0)}W
 \le C\|W\|_{L^2(B_{4\rho}(x_0))}\le C.
\]
Here the constant depends only on \(n,\ell,\rho\) and the displayed
bound for \(B\); see \cite{Trudinger}. In particular
\(\|W\|_\infty\le C\) uniformly on \(\Omega\).
The right-hand sides of \eqref{eq:pen} are consequently locally
bounded. Interior \(W^{2,l}\) estimates for every finite \(l\)
give \(C^{1,\theta}_{\rm loc}\) regularity for \(0<\theta<1\).
On bounded sets the powers are H\"older continuous on the
nonnegative quadrant, while \(h_\delta\) is locally Lipschitz in
\((y,s)\). Schauder estimates therefore give \(u,v\in C^2(\Omega)\).

The Harnack inequality on the same balls now reads
\[
 \sup_{B_\rho(x_0)}W\le C_H\inf_{B_\rho(x_0)}W,
\]
with \(C_H\) independent of \(x_0\) and \(\delta\). Since
\(G(U)>0\), \(W\) is not identically zero, and propagation along
overlapping balls on the connected cylinder gives \(W>0\).
At any fixed \(y\), the circle \(\{y\}\times\T\) is covered by a
chain of at most \(N_\rho\) overlapping balls of radius \(\rho\),
where \(N_\rho\) depends only on \(\rho\) and the period \(2\pi\).
Iteration of Harnack along the chain yields
\begin{equation}\label{eq:harnack}
 \sup_{t\in\T} W(y,t)\le C\inf_{t\in\T}W(y,t)\le C M(y),
 \qquad M=\avg W .
\end{equation}
Since \(M=2w\), \eqref{eq:radial} shows that \(W\to0\) uniformly
as \(|y|\to\infty\), uniformly also in \(\delta\).

\subsection*{The antisymmetric mode}
Subtracting the equations gives \(-\Delta D=A D\), where
\begin{equation}\label{eq:A}
 0\le A\le C W^{p-1}.
\end{equation}
Indeed the coefficient is \(\lambda_\delta\) times the sum of
\((u^q+v^q)(u^a-v^a)/(u-v)\) and the divided difference of \(h_\delta\).
The first quotient involving \(u^a-v^a\) is bounded by \(aW^q\).
On \([0,W]\) the truncated function is increasing with Lipschitz
constant at most \(p\delta W^{p-1}\). This also bounds the slope
on the linear branch, which is active only when
\(e\omega(y)\le\delta s^{p-1}\). Where \(D=0\) one may set \(A=0\).
There is a fixed \(R_0\ge1\) such that \(A\le1/2\) for \(|y|\ge R_0\).

Now \(D(y,t+\pi)=-D(y,t)\), so
\(\avg{|\partial_tD|^2}\ge\avg{D^2}\). For \(K(y)=\avg{D^2}(y)\),
periodic integration by parts gives on \(|y|\ge R_0\)
\[
 \tfrac12\Delta_y K
 =\avg{|\nabla_yD|^2+|\partial_tD|^2-A D^2}\ge\tfrac12 K .
\]
The function \(K\) tends to zero by \eqref{eq:harnack}.
For \(Z(y)=C\exp(-(|y|-R_0)/2)\) one has
\[
 (-\Delta_y+1)Z
 =\left(\frac34+\frac{k-1}{2|y|}\right)Z\ge0.
\]
Choose \(C\) to dominate \(K\) on \(|y|=R_0\), uniformly in
\(\delta\). The function \(K-Z\) has nonpositive boundary values,
tends to zero at infinity, and satisfies
\((-\Delta_y+1)(K-Z)\le0\). If its supremum were positive,
it would be attained at an interior point, contradicting this
inequality. Therefore
\begin{equation}\label{eq:K}
 K(y)\le C e^{-c|y|}\quad (|y|\ge R_0).
\end{equation}

\subsection*{The mean and the original coefficients}
By \eqref{eq:gap}, \(0<q<1\). Thus
\[
 |u^q-v^q|\le |D|^q,\qquad |u^a-v^a|\le a W^q |D|\le C|D|.
\]
Average \eqref{eq:sum} and use
\(\avg{|D|^{q+1}}\le K^{(q+1)/2}\). Equations \eqref{eq:lambda}
and \eqref{eq:K} imply
\begin{equation}\label{eq:mean}
 -\Delta_y M\le C_0\omega(y) M+C_1 e^{-c_1|y|}
 \qquad (|y|\ge R_0)
\end{equation}
with uniform constants. Put \(d_\gamma=\gamma(k-2-\gamma)>0\).
Increase \(R_0\), independently of \(\delta\), so that
\(C_0\omega(y)\le\min(d_\gamma,H_k)/(2|y|^2)\) on the exterior.
Since \(-\Delta_y |y|^{-\gamma}=d_\gamma |y|^{-\gamma-2}\),
one can choose \(C_2\) so that \(P(y)=C_2|y|^{-\gamma}\) dominates
\(M\) on the inner sphere and
\((-\Delta_y-C_0\omega)P\ge C_1 e^{-c_1|y|}\) throughout the exterior.
We justify the comparison in the homogeneous energy space.
All unlabelled integrals in this comparison are over \(\R^k\).
Put \(\mathcal E=\{|y|>R_0\}\), \(V=C_0\omega\), and \(Q=M-P\).
Subtracting the two differential inequalities gives
\[
 -\Delta_y Q\le VQ\quad\hbox{in }\mathcal E,\qquad
 Q\le0\quad\hbox{on }\partial\mathcal E,\qquad
 0\le V\le H_k/(2|y|^2).
\]
The mean \(M\) belongs to \(\dot H^1(\R^k)\).
The condition \(\gamma>(k-2)/2\) gives
\[
 \int_{\mathcal E}|\nabla P|^2
 =C\int_{R_0}^\infty r^{k-3-2\gamma}\,dr<\infty,
 \qquad
 \int_{\mathcal E}|P|^{s_k}
 =C\int_{R_0}^\infty r^{k-1-\gamma s_k}\,dr<\infty.
\]
Thus \(\phi=Q_+\), extended by zero to \(B_{R_0}\), has zero
trace on the inner sphere, finite Dirichlet energy, and belongs to
\(L^{s_k}(\R^k)\).

For \(R>2R_0\), choose a smooth cutoff \(\chi_R\) equal to one
on \(B_R\), zero outside \(B_{2R}\), with
\(0\le\chi_R\le1\) and \(|\nabla\chi_R|\le C/R\).
The function \(\chi_R^2\phi\) is in \(H^1_0\) of a bounded
exterior annulus. Approximation there by nonnegative smooth
compactly supported functions makes it a valid test in the
difference inequality. Using
\(\nabla Q\cdot\nabla\phi=|\nabla\phi|^2\) and
\(Q\phi=\phi^2\) almost everywhere gives
\[
 \int|\nabla(\chi_R\phi)|^2
 \le\int V(\chi_R\phi)^2+\int\phi^2|\nabla\chi_R|^2
 \le\frac12\int|\nabla(\chi_R\phi)|^2
       +\int\phi^2|\nabla\chi_R|^2,
\]
where Hardy's inequality applies to the zero extension of
\(\chi_R\phi\). The remaining error satisfies
\[
 \int\phi^2|\nabla\chi_R|^2
 \le CR^{-2}\int_{R<|y|<2R}\phi^2
 \le C\left(\int_{R<|y|<2R}|\phi|^{s_k}\right)^{2/s_k}
 \longrightarrow0.
\]
In the second inequality, the volume factor is \(O(R^2)\)
because \(1-2/s_k=2/k\). Letting \(R\to\infty\) first on
each fixed ball shows that \(\nabla\phi=0\). The
\(L^{s_k}\) integrability then implies \(\phi=0\).
This proves \(M\le P\). The same cutoff estimate, together with
the vanishing tail of \(\int|\nabla\phi|^2\), also justifies
approximation in the Dirichlet norm.
Returning to \eqref{eq:harnack} proves
\[
 W(y,t)\le C(1+|y|)^{-\gamma}\qquad\text{throughout }\Omega.
\]
The strict inequality \(\gamma(p-1)>2+\sigma\) in \eqref{eq:choice}
is what removes the truncation. In fact
\(\sup_{\Omega}W^{p-1}/\omega\le C_*\) independently of \(\delta\).
For every \(0<\delta<\delta_0:=\min(1,e/(1+C_*))\), this implies
\(h_\delta(y,u)=\delta u^p\) and \(h_\delta(y,v)=\delta v^p\) everywhere.

The two coefficients in \eqref{eq:pen} are consequently
\(\lambda_\delta(1+\delta)\) on the pure powers and
\(-\lambda_\delta\) on the mixed powers. Multiplication of both
components by \([\lambda_\delta(1+\delta)]^{1/(p-1)}\) gives
\(\mu_1=\mu_2=1\) and
\[
 \beta=-\frac1{1+\delta}.
\]
This realizes every \(-1<\beta<-1+\eta\) with
\(\eta=\delta_0/(1+\delta_0)\). Periodic extension gives a classical
solution on \(\R^n\). The shift and \(G>0\) ensure that both components
are nonzero and unequal. The cell energy is finite and positive, so
its sum over all translated cells diverges. This proves
the normalized case of Theorem~\ref{thm:exist}.

Finally, let \(\mu>0\) and
\(-\mu<\beta<-\mu+\mu\eta\). The normalized construction applies
to \(b=\beta/\mu\). If \((u_0,v_0)\) is its solution with
self-coefficients one and coupling \(b\), define
\[
 u=\mu^{-1/(p-1)}u_0,\qquad
 v=\mu^{-1/(p-1)}v_0.
\]
Since both the pure and mixed terms have degree \(p\),
\[
 -\Delta u=\mu u^p+\mu b\,u^qv^{q+1}
           =\mu u^p+\beta u^qv^{q+1},
\]
and the same computation gives the second equation.
All the stated symmetry, decay, and energy properties are
preserved by this positive constant multiplication. This completes
the proof of Theorem~\ref{thm:exist}.
\qed 
\begin{remark}[The endpoint for periodic solutions with general coefficients]
The restriction \(p>p_*(n)\) is necessary in the one-periodic class
for \emph{arbitrary} \(\mu_1,\mu_2>0\), provided
\(-\sqrt{\mu_1\mu_2}<\beta<0\). Indeed, set
\[
 U=\mu_1^{1/(2a)}u,\quad V=\mu_2^{1/(2a)}v,\quad
 d_i=\mu_i^{-1/a},\quad \theta=-\frac{\beta}{\sqrt{\mu_1\mu_2}}.
\]
For \(-\sqrt{\mu_1\mu_2}<\beta<0\), the normalized equations are
\begin{equation}\label{eq:unequal-normalization}
 -d_1\Delta U=U^p-\theta U^qV^a,\qquad
 -d_2\Delta V=V^p-\theta V^qU^a.
\end{equation}
Put \(S=d_1U+d_2V\).
Since \(0<\theta<1\), addition gives
\[
 \begin{split}
 -\Delta S
 &= (1-\theta)(U^p+V^p)
    +\theta(U^q-V^q)(U^a-V^a)\\
 &\ge (1-\theta)(U^p+V^p)\ge cS^p.
 \end{split}
\]
A nontrivial nonnegative pair therefore has \(S>0\) by the
strong minimum principle. Periodic integration and Jensen's
inequality show that \(M(y)=\avg S(y)>0\) satisfies
\(-\Delta_yM\ge cM^p\) on \(\R^k\), \(k=n-1\).
This argument needs neither transverse radial symmetry nor
boundedness of the entire solution.

For completeness, let \(b(r)\) be the spherical mean of \(M\)
about the origin. Jensen's inequality gives
\[
 -(r^{k-1}b')'\ge cr^{k-1}b^p,\qquad b'(0)=0.
\]
Thus \(b'\le0\). For \(r\ge1\), integration on \([r/2,r]\)
and monotonicity give
\(-b'(r)\ge c_1 r b(r)^p\). Integration of
\((b^{1-p})'\ge(p-1)c_1r\) yields
\(b(r)\le C r^{-2/(p-1)}\) for \(r\ge2\).
On the other hand, retaining the positive integral over
\([0,1]\) gives \(-r^{k-1}b'(r)\ge c_2>0\) for \(r\ge1\).
Consequently
\[
 b(r)\ge b(r)-b(2r)
       \ge c_2\int_r^{2r}s^{1-k}\,ds
       \ge c_3 r^{2-k}.
\]
The two bounds are incompatible when \(1<p<k/(k-2)\).
At \(p=k/(k-2)\), the lower bound instead implies
\[
 -r^{k-1}b'(r)\ge c_4\int_1^r s^{k-1-(k-2)p}\,ds
                    =c_4\log r.
\]
Integrating on \([r,2r]\) now gives
\(b(r)\ge c_5r^{2-k}\log r\), contradicting the upper bound
\(Cr^{2-k}\). Hence \(p>k/(k-2)=p_*(n)\) is necessary for a
nontrivial one-periodic solution. The endpoint question for
entire solutions without this periodicity assumption remains
outside this argument.
\end{remark}

\section{Weak competition and the scalar limit}
\begin{proof}[Proof of Theorem~\ref{thm:weak}]
The case \(\beta=0\) is the scalar theorem \cite{GS} applied to each
component. If the assertion failed, there would be negative
\(\beta_j\to0\), exponents \(q_j\in(0,\bar q]\), and nontrivial
nonnegative entire solutions \((u_j,v_j)\) for these parameters.
After extraction, \(q_j\to q_\infty\in[0,\bar q]\).
Put \(S_j=u_j+v_j\). We select points so that the normalized
solutions have a bound independent of \(q_j\).
Start at any point where \(S_j>0\). At a current point \(x\), stop if
\[
 S_j(y)\le2S_j(x)
 \quad\text{for all }|y-x|\le jS_j(x)^{-q_j}.
\]
Otherwise move to a point of this ball where \(S_j\) more than
doubles. If \(x^{(0)}\) is the starting point, an infinite sequence
of moves would have total length at most
\[
 jS_j(x^{(0)})^{-q_j}\sum_{h=0}^\infty 2^{-h q_j}
 =\frac{jS_j(x^{(0)})^{-q_j}}{1-2^{-q_j}}<\infty.
\]
All selected points would then lie in a compact ball, whereas
\(S_j\) would tend to infinity along them. Continuity rules this
out. The argument uses \(q_j>0\) for each \(j\), and requires no
uniform bound for the displayed total length. It is a version of
the doubling selection in \cite[Lemma~5.1]{PQS}.
Denote the stopping point by \(x_j\), and set
\[
 A_j=S_j(x_j)>0,\qquad L_j=A_j^{-q_j}.
\]

The pair
\[
 (U_j,V_j)(z)=A_j^{-1}(u_j,v_j)(x_j+L_j z)
\]
solves the same system with exponent \(q_j\) and coefficient
\(\beta_j\), because \(L_j^2A_j^{2q_j}=1\). It satisfies
\[
 (U_j+V_j)(0)=1,\qquad 0\le U_j+V_j\le2
 \quad\text{on }B_j(0).
\]
All nonlinearities on the right-hand sides are uniformly bounded
on each fixed ball, also when \(q_j\to0\). For example,
\[
 |\beta_j|U_j^{q_j}V_j^{q_j+1}
 \le|\beta_j|\,2^{2\bar q+1}\longrightarrow0.
\]
Interior \(W^{2,l}\) estimates with \(l>n\), compact embedding and
diagonal extraction give local \(C^1\) convergence along a subsequence.
The functions \(s\mapsto s^{1+2q_j}\) converge uniformly on
\([0,2]\) to \(s\mapsto s^{1+2q_\infty}\), including when
\(q_\infty=0\). Passing to the limit in distributions gives
\[
 -\Delta U=\mu_1 U^{1+2q_\infty},\qquad
 -\Delta V=\mu_2 V^{1+2q_\infty},\qquad
 (U+V)(0)=1.
\]
Both limits are nonnegative and classical by interior regularity.
If \(q_\infty>0\), then \(1<1+2q_\infty<p_S(n)\).
Multiplying the \(i\)-th component by \(\mu_i^{1/(2q_\infty)}\)
reduces its equation to the scalar equation in \cite{GS}.
Thus \(U=V=0\), contrary to the normalization.

If \(q_\infty=0\), the equations become
\(-\Delta U=\mu_1U\) and \(-\Delta V=\mu_2V\).
For completeness, no nonzero nonnegative entire solution \(Z\)
of \(-\Delta Z=\mu Z\), \(\mu>0\), can exist.
Choose \(R\) so large that the first Dirichlet eigenvalue
\(\lambda_1(B_R)\) is smaller than \(\mu\), and let
\(\varphi_R>0\) be its eigenfunction. Green's formula gives
\[
 (\mu-\lambda_1(B_R))\int_{B_R}Z\varphi_R
 =\int_{\partial B_R}Z\,\partial_\nu\varphi_R\le0.
\]
Here \(\varphi_R=0\) and \(\partial_\nu\varphi_R<0\) on the
boundary. Hence \(Z=0\) in \(B_R\); letting \(R\) increase gives
\(Z=0\) everywhere. Applying this to \(U\) and \(V\) again
contradicts \((U+V)(0)=1\).
Shrinking \(\varepsilon_*\) if necessary keeps the interval strictly
above \(-\sqrt{\mu_1\mu_2}\).
\end{proof}

\begin{remark}[Dependence on the exponent]\label{rem:uniform}
Theorem~\ref{thm:weak} provides a constant independent of the
individual exponent \(q\) throughout \(0<q\le\bar q\).
In dimensions \(n=3,4\), the stronger nonexistence result of
\cite{QS12}, recalled in the introduction, permits the choice
\(\varepsilon_*=\tfrac12\sqrt{\mu_1\mu_2}\) for the entire interval
\(0<q<2/(n-2)\).
In dimensions \(n\ge5\), uniformity up to this upper endpoint is
not established here. If \(q_j\to2/(n-2)\), the normalized limit
in the preceding proof solves the critical scalar equations, which
admit, for example,
\[
 U(x)=\left(1+\frac{\mu_1|x|^2}{n(n-2)}\right)^{-(n-2)/2},
 \qquad V(x)=0.
\]
Direct differentiation gives
\(-\Delta U=\mu_1U^{(n+2)/(n-2)}\), and this pair satisfies
\((U+V)(0)=1\) and \(0\le U+V\le1\).
It is therefore compatible with all the normalization and
compactness conclusions above. The possible occurrence of this
limit prevents that argument from proving a constant uniform
over the full subcritical interval. 
\end{remark}

\medskip 

\noindent\textbf{Acknowledgements}  The research of H. Yang is supported by NSFC 12671139 and 12301140. The authors used AI-assisted tools, principally ChatGPT/Codex, for brainstorming and manuscript editing. AI output was not treated as a mathematical authority. The authors take full responsibility for all mathematical statements, arguments, citations, and conclusions in this paper.

\end{document}